\documentclass{gITR2e}
\usepackage{eucal}

\theoremstyle{plain}
\newtheorem{theorem}{Theorem}[section]

\newtheorem{proposition}[theorem]{Proposition}

\makeatletter
\def\ps@plain{%
     \let\@mkboth\@gobbletwo
     \let\@oddfoot\@evenfoot
     \def\@oddhead{\hbox to \textwidth{}}%
     \let\@evenhead\@oddhead
}
\makeatother

\begin{document}

\title{{\itshape Integrals of Lipschitz--Hankel type, Legendre functions and\\ table errata}}
\author{\name{Robert S. Maier$^{\ast}$\thanks{$^\ast$Email: rsm@math.arizona.edu}}
\affil{Depts.\ of Mathematics and Physics, University of Arizona, Tucson, AZ 85721, USA}}

\def\draftnote{}%

\maketitle

\begin{abstract}
  The complete Lipschitz--Hankel integrals (LHIs) include the Laplace
  transforms of the Bessel functions, multiplied by powers.  Such Laplace
  transforms can be evaluated using associated Legendre functions.  It is
  noted that there are errors in published versions of these evaluations,
  and a merged and emended list of seven transforms is given.  Errata for
  standard reference works, such as the table of Gradshteyn and Ryzhik, are
  also given.  Most of the errors are attributable to inconsistent
  normalization of the Legendre functions.  These transforms can be viewed
  as limits of incomplete LHIs, which find application in communication
  theory.
\end{abstract}


\begin{classcode}
33C10; 44A10  
\end{classcode}

\section{Introduction}

A complete Lipschitz--Hankel integral (LHI) is an integral transform of a
product of one or more Bessel or cylindrical functions.  The simplest such
integrals are the Laplace transforms $\mathcal{L}\{t^\nu {\mathcal
  C}_\mu(t)\}(s)$, where ${\mathcal C}_\mu = J_\mu,Y_\mu,I_\mu,K_\mu$
specifies a Bessel or modified Bessel function, of order~$\mu$.  It has
long been known that each of these Laplace transforms can be expressed
in~terms of an (associated) Legendre function~\cite[\S\,13.21]{Watson44}.

In published tables \cite{Erdelyi54,Roberts66,Gradshteyn2015} and other
reference works \cite{Watson44,Agrest71}, erroneous evaluations of these
integrals have been found; these are reported in~\S\,\ref{sec:errata}.
Of~the errors in
\cite{Erdelyi54,Roberts66,Gradshteyn2015,Watson44,Agrest71}, most are
attributable to inconsistent normalization of the Legendre functions; in
particular, of the second Legendre function~$Q_\nu^\mu(z)$.  None of
\cite{Erdelyi54,Roberts66,Gradshteyn2015,Watson44,Agrest71}, and at least
three other tables in which these integrals appear correctly
\cite{Magnus66,Prudnikov86,Olver2010}, includes a complete set of
evaluations.  But a complete set of seven can be assembled by merging
transforms from different sources, emending as necessary.  Such a set is
supplied in Theorems \ref{thm:1} and~\ref{thm:2} below.

These theorems may prove useful in verifying the entries in databases of
definite integrals, and more generally in avoiding errors in symbolic
computation.  Of~the seven transforms, (\ref{eq:J})~and~(\ref{eq:Ip}) have
recently been used in geometric analysis~\cite[\S\,8.8]{Taylor2011b}, and
(\ref{eq:Ip})~has been used in stochastic analysis~\cite{Matsumoto2002} and
path integration~\cite{Grosche92}.  The others also have applications.
Each is given here in a compact, trigonometrically parametrized form.

Incomplete LHIs (ILHIs), such as the incomplete integrals
$\mathcal{C}e_{\nu,\mu}(T;s)$ that have the Laplace transforms
$\mathcal{L}\{t^\nu {\mathcal C}_\mu(t)\}(s)$ as their $T\to\infty$ limits,
have recently been applied in electromagnetics~\cite{Dvorak94} and
communication theory~\cite{Sofotasios2014}. In particular,
$Je_{\nu,\mu}(T;s)$ and $Ie_{\nu,\mu}(T;s)$ have been applied.  For
comparison purposes, explicit expressions for $\mathcal{L}\{t^\nu
J_\mu(t)\}(s)$ and $\mathcal{L}\{t^\nu I_\mu(t)\}(s)$ in~terms of Legendre
functions, not trigonometrically parametrized, are given
in~\S\,\ref{sec:suppl}.  The transform $\mathcal{L}\{t^\nu K_\mu(t)\}(s)$
is also commented~on.

\section{Legendre functions as Laplace transforms}

The Legendre functions of degree~$\nu$ and order~$\mu$, which are solutions
of the corresponding Legendre differential equation, include $P_\nu^\mu(z)$
and $Q_\nu^\mu(z)$, which are defined on the complement of the cut
$z\in(-\infty,1]$.  They also include the Ferrers functions
  ${\rm{P}}_\nu^\mu(z)$ and ${\rm Q}_\nu^\mu(z)$, sometimes called Legendre
  functions `on~the cut,' which are defined on the complement of
  $z\in(-\infty,-1]\cup[1,\infty)$ and occur most often with $z\in(-1,1)$;
      especially, in physical applications.  In the absence of cuts, all
      these functions would be multivalued.

In the older English-language literature (e.g.,~\cite{Hobson31}),
$P_\nu^\mu,Q_\nu^\mu$ and ${\rm P}_\nu^\mu,{\rm Q}_\nu^\mu$ were not
typographically distinguished, and this remains true in a few recent works
(e.g.,~\cite{Gradshteyn2015}).  In the German literature
(e.g.,~\cite{Lense54,Magnus66}), they are denoted by ${\mathfrak
  P}_\nu^\mu,{\mathfrak Q}_\nu^\mu$ and ${P}_\nu^\mu,{Q}_\nu^\mu$.  In the
recent NIST handbook~\cite{Olver2010}, ${\rm P}_\nu^\mu,{\rm Q}_\nu^\mu$
are denoted by $\mathsf{P}_\nu^\mu,\mathsf{Q}_\nu^\mu$, to stress the
dissimilarity to~${P}_\nu^\mu,{Q}_\nu^\mu$.

The function $Q_\nu^\mu(z)$, defined on the complement of the cut
$z\in(-\infty,1]$, has historically been defined in three ways, differing
    only in normalization.  The most recent (and perhaps the best)
    definition is that of the late FWJ Olver~\cite{Olver97}.  But for many
    decades, the standard definition in the literature has been that of
    EW~Hobson~\cite{Hobson31}.  There is also a definition originating with
    EW~Barnes~\cite{Barnes08}, which is no~longer used but was employed in
    at~least two major works that are still consulted
    \cite{Bateman32,Watson44}.  The three definitions are related by
    \begin{displaymath}
      e^{-\mu\pi{\rm i}}\,[Q_\nu^\mu]_H = \Gamma(\nu+\mu+1)\,[Q_\nu^\mu]_O
      = \frac{\sin(\nu\pi)}{\sin[(\nu+\mu)\pi]}\,[Q_\nu^\mu]_B,
    \end{displaymath}
    the subscript indicating the originator.  In the NIST
    handbook~\cite{Olver2010}, the function of Olver is denoted
    by~${\boldsymbol{Q}}_\nu^\mu$, the symbol~${Q}_\nu^\mu$ being reserved
    for the standard (Hobson) function.  But the Barnes function, when in
    use, was never distinguished typographically from that of Hobson.  It
    should be noted that alternative definitions of $P_\nu^\mu$ as~well
    have been employed (e.g., in~\cite{Snow52,Frolov2001}).  However,
    $Q_\nu^\mu$~is the chief potential source of confusion.

The following are the complete set of Laplace transforms, employing
Hobson's now-standard definition of~$Q_\nu^\mu$.  They are parametrized to
agree with the circular- and hyperbolic-trigonometric versions given
in~\cite{Watson44}.  (Except for~\cite{Gradshteyn2015}, the other cited
works give unnormalized, non-trigonometric forms, but their forms easily
reduce to the following.)  Each identity has been checked by numerical
quadrature.

\begin{theorem}
\label{thm:1}
The following identities hold, when\/ $\theta\in(0,\pi/2)$ for (\ref{eq:J})
and~(\ref{eq:Y}), resp.\ $\xi\in(0,\infty)$ for (\ref{eq:I})
and~(\ref{eq:K}), and\/ $\xi\in\mathbb{R}$ for~(\ref{eq:K2}).  \begingroup
\allowdisplaybreaks
\begin{align}
\label{eq:J}
  \int_0^\infty e^{-t\cos\theta}\:t^\nu J_\mu(t\sin\theta)\,{\rm d}t
& = \hphantom{(-2/\pi)}\Gamma(\nu+\mu+1)\,{\rm P}_\nu^{-\mu}(\cos\theta),\\
\label{eq:Y}
  \int_0^\infty e^{-t\cos\theta}\:t^\nu\, Y_\mu(t\sin\theta)\,{\rm d}t
& = {-(2/\pi)}\,\Gamma(\nu+\mu+1)\,{\rm
    Q}_\nu^{-\mu}(\cos\theta), \\
\label{eq:I}
  \int_0^\infty e^{-t\cosh\xi}\:t^\nu I_\mu(t\sinh\xi)\,{\rm d}t
& = \hphantom{(-2/\pi)}\Gamma(\nu+\mu+1)\,{P}_\nu^{-\mu}(\cosh\xi),\\
\label{eq:K}
  \int_0^\infty e^{-t\cosh\xi}\:t^\nu K_\mu(t\sinh\xi)\,{\rm d}t
& = \hphantom{(-2/\pi)}\Gamma(\nu+\mu+1)\,e^{\mu\pi{\rm i}}\,{Q}_\nu^{-\mu}(\cosh\xi),\\
\label{eq:K2}
  \int_0^\infty e^{-t\sinh\xi}\:t^\nu K_\mu(t\cosh\xi)\,{\rm d}t & =\nonumber\\
{(\pi/2)^{1/2}}\,\Gamma(\nu+\mu+1)&\Gamma(\nu-\mu+1)\,(\cosh\xi)^{-1/2}\,{\rm P}_{\mu-1/2}^{-\nu-1/2}(\tanh\xi).
\end{align}%
\endgroup
It is assumed that the integrals converge, for which one needs ${\rm
  Re}(\nu+\nobreak\mu)>\nobreak-1$ in (\ref{eq:J}) and~(\ref{eq:I}), and\/ 
${\rm Re}(\nu\pm\mu)>\nobreak-1$ in (\ref{eq:Y}), (\ref{eq:K})
and~(\ref{eq:K2}).
\end{theorem}

The five identities of the theorem can be extended to the complex domain
(see the cited references).  Also, (\ref{eq:I}),(\ref{eq:K}) are equivalent
to the following.
\begin{theorem}
\label{thm:2}
\begin{align}
  \label{eq:Ip}
  \int_0^\infty e^{-t\cosh\xi}\:t^\nu I_\mu(t\sinh\xi)\,{\rm d}t&=
   (2/\pi)^{1/2}(\sinh\xi)^{-1/2}\,e^{-(\nu+1/2)\pi{\rm i}}\,{Q}_{\mu-1/2}^{\nu+1/2}(\coth\xi),\\
  \int_0^\infty e^{-t\cosh\xi}\:t^\nu K_\mu(t\sinh\xi)\,{\rm d}t&=\nonumber\\
  (\pi/2)^{1/2}\, \Gamma(\nu+\mu+1)&\Gamma(\nu-\mu+1)\,
  (\sinh\xi)^{-1/2}\,{P}_{\mu-1/2}^{-\nu-1/2}(\coth\xi).
\label{eq:Kp}
\end{align}
\end{theorem}
These two supplementary identities are obtained by rewriting the right
sides of (\ref{eq:I}),(\ref{eq:K}) with the aid of Whipple's formula
\cite[Ex.~12.2]{Olver97} and the fact that ${\boldsymbol Q}_\nu^{\pm\mu}$,
i.e., $[Q_\nu^{\pm \mu}]_O$, are equal to each other.

Of the seven Laplace transform evaluations listed in Theorems \ref{thm:1}
and~\ref{thm:2}, the number appearing, in whatever form, in
\cite{Watson44},\cite{Erdelyi54},\cite{Roberts66},\cite{Magnus66},\cite{Agrest71},\cite{Oberhettinger73},\cite{Prudnikov86},\cite{Olver2010},\cite{Gradshteyn2015}
is respectively 6,6,6,4,4,7,5,4,6.  The number appearing \emph{incorrectly}
is respectively 1,4,2,0,1,0,0,0,3.  (The ordering used here is
chronological.)  In some cases the right sides are given as expressions
involving the Gauss hypergeometric function~${}_2F_1$, but can be rewritten
using Legendre functions so as to agree with entries in the foregoing list.

The identities (\ref{eq:J})--(\ref{eq:K}), in particular (\ref{eq:J})
and~(\ref{eq:I}), are proved by series expansion and term-by-term
integration; which is a procedure that can be traced to Hankel's 1875 proof
of~(\ref{eq:J}), and extended to more general integral transforms (see,
e.g.,~\cite{Srivastava79}).  The most troublesome of them may
be~(\ref{eq:Y}).  Only in~\cite{Agrest71} is the Laplace transform of
$t^\nu Y_\mu(t)$ given in the form based on the second Ferrers
function~${\rm Q}_\nu^{-\mu}$ that appears in~(\ref{eq:Y}).  In many works
this transform is given as a combination of the Ferrers functions ${\rm
  P}_\nu^{\mu},{\rm P}_\nu^{-\mu}$, but by standard Legendre identities
\cite{Olver97,Olver2010} this can be seen to be equivalent to~(\ref{eq:Y}).
The little-known evaluation~(\ref{eq:K2}) appears only in
\cite{Magnus66,Oberhettinger73}, where it is given correctly, in a
non-trigonometric parametrization.  (See p.~92, fourth equation down;
resp.\ p.~151, eq.~1.15.16.)

It is evident from (\ref{eq:J})--(\ref{eq:K}) that the interchange symmetry
between the pairs $P_\nu^{-\mu},Q_\nu^{-\mu}$ and ${\rm P}_\nu^{-\mu},{\rm
  Q_\nu^{-\mu}}$ is not reflected in a symmetry between the pairs
$J_\mu,Y_\mu$ and~$I_\mu,K_\mu$.  The difficulty is the prefactor $2/\pi$
in~(\ref{eq:Y}): there is no such prefactor in~(\ref{eq:K}).  This
asymmetry is attributable to the definition of the second modified Bessel
function~$K_\mu$, `MacDonald's function,' which includes a factor equal
to~$\pi/2$.  It has been suggested that using a differently normalized
function $Kh_\mu$, equal to $2/\pi$ times~$K_\mu$, would be more
reasonable~\cite{Jeffreys56}.  (The~`h' is for Heaviside, who preferred
this normalization.)  Rewriting (\ref{eq:K}),(\ref{eq:K2})
and~(\ref{eq:Kp}) in~terms of~$Kh_\mu$ would exhibit a symmetry between
$J_\mu,Y_\mu$ and~$I_\mu,Kh_\mu$.

\section{Some consequences}
\label{sec:suppl}

The incomplete Lipschitz--Hankel integral $\mathcal{C}e_{\nu,\mu}(T;s)$,
for $\mathcal{C}=J,Y,I,K$, is defined for real~$s$, when $T>0$, by
\begin{equation}
\mathcal{C}e_{\nu,\mu}(T;s) = \int_0^T e^{-st} \,t^\nu
\mathcal{C}_\mu(t)\,dt, 
\end{equation}
provided that the integral converges.  For convergence, one needs ${\rm
  Re}(\nu+\mu)>-1$ for $\mathcal{C}=J,Y$, and ${\rm Re}(\nu\pm\mu)>-1$ for
$\mathcal{C}=I,K$.  The limit as $T\to\infty$ of
$\mathcal{C}e_{\nu,\mu}(T;s)$ is the Laplace transform $\mathcal{L}\{t^\nu
{\mathcal C}_\mu(t)\}(s)$, provided that this transform exists.  The
transform exists if $s>0$ for $\mathcal{C}=J,Y$, if $s>1$ for
$\mathcal{C}=I$, and if $s>-1$ for $\mathcal{C}=K$.

For comparison with previous work on ILHIs
\cite{Agrest71,Dvorak94,Sofotasios2014}, the following expressions for
$\mathcal{L}\{t^\nu {\mathcal C}_\mu(t)\}(s)$ with $\mathcal{C}=J,I,K$ may
be useful.  They are non-trigonometric versions of
(\ref{eq:J}),(\ref{eq:I}),(\ref{eq:K2}),(\ref{eq:Kp}), respectively, and
none involves the troublesome functions $Q_\nu^\mu,{\rm Q}_\nu^\mu$.

\begin{align}
  \mathcal{L}\{t^\nu {J}_\mu(t)\}(s)&= \Gamma(\nu+\mu+1)\: (1+s^2)^{-(\nu+1)/2}
          \:{\rm P}_\nu^{-\mu}(s/\sqrt{1+s^2}),\qquad s\in(0,\infty);\label{eq:JJ}\\
  \mathcal{L}\{t^\nu {I}_\mu(t)\}(s)&= \Gamma(\nu+\mu+1)\: (s^2-1)^{-(\nu+1)/2}
          \:{P}_\nu^{-\mu}(s/\sqrt{s^2-1}),\qquad s\in(1,\infty);\label{eq:II}\\
  \mathcal{L}\{t^\nu {K}_\mu(t)\}(s)&= (\pi/2)^{1/2}\,\Gamma(\nu+\mu+1)\Gamma(\nu-\mu+1)\nonumber\\
   &\quad\qquad\qquad\qquad {}\times\left\{
\begin{aligned}
  (1-s^2)^{-(2\nu+1)/4}\:{\rm P}_{\mu-1/2}^{-\nu-1/2}(s), & \qquad s\in(-1,1),\label{eq:KK}\\
  (s^2-1)^{-(2\nu+1)/4}\:P_{\mu-1/2}^{-\nu-1/2}(s), & \qquad s\in(1,\infty).
\end{aligned}
\right.
\end{align}
The conditions on $\nu,\mu$ for convergence are as stated.

The standard reference on ILHIs~\cite{Agrest71} includes versions of
(\ref{eq:JJ}) and~(\ref{eq:II}).  But as is noted in~\S\,\ref{sec:errata}
below, its version of~(\ref{eq:II}) includes an erroneous $-(2/\pi)$
factor.  Also, the Laplace transform of $t^\nu K_\mu(t)$ given
in~\cite{Agrest71} is valid only when $s\in(1,\infty)$; and involving
$Q_\nu^{-\mu}$ rather than $P_{\mu-1/2}^{-\nu-1/2}$, it is a
non-trigonometric version not of~(\ref{eq:Kp}), but rather of~(\ref{eq:K}).

The duality between the expressions for $\mathcal{L}\{t^\nu
{K}_\mu(t)\}(s)$ on $s\in(-1,1)$ and $s\in(1,\infty)$, evident
in~(\ref{eq:KK}), deserves to be better known.  The value at~$s=1$ must be
supplied separately.
\begin{proposition}
\label{prop:31}
At\/ $s=1$, the Laplace transform of\/ $t^\nu K_\mu(t)$ is given by
\begin{equation}
\label{eq:s1}
\mathcal{L}\{t^\nu {K}_\mu(t)\}(s=1) = 
\frac{\pi^{1/2}\, \Gamma(\nu+\mu+1)\Gamma(\nu-\mu+1)}{2^{\nu+1}\,\Gamma(\nu+3/2)},
\end{equation}
where one needs\/ ${\rm Re}(\nu\pm\mu)>-1$ for convergence.
\end{proposition}
This comes from~(\ref{eq:KK}) by considering the asymptotic behavior of
${\rm P}_{\mu-1/2}^{-\nu-1/2}(s)$ as~$s\to1^-$, or of
${P}_{\mu-1/2}^{-\nu-1/2}(s)$ as~$s\to1^+$, which are well
known~\cite{Olver97}.  Viewed as a function of~$\nu$, the right side
of~(\ref{eq:s1}) is the Mellin transform of~$e^{-t}K_\mu(t)$, and it has
been tabulated as~such \cite[eq.~6.8(28)]{Erdelyi54}.

\section{Table errata}
\label{sec:errata}

The following are the errata.  None seems to have been previously reported,
though the incompatibility between the definitions used in
\cite{Watson44,Gradshteyn2015} has been noted~\cite{Matsumoto2002}.

\bigskip
\noindent
\textsc{G. N. Watson}, \emph{A Treatise on the Theory of Bessel Functions},
2nd ed., Cambridge University Press, Cambridge, UK, 1944.  Cited here as
Ref.~\cite{Watson44}.

\nopagebreak
\medskip
The Barnes definition of~$Q_\nu^\mu$ is used, and ${\rm P}_\nu^\mu,{\rm
  Q}_\nu^\mu$ are not visually distinguished from
${P}_\nu^\mu,{Q}_\nu^\mu$.  On pp.\ 387--388, in~\S\,13.21, versions of
identities
(\ref{eq:J}),(\ref{eq:Y}),\allowbreak(\ref{eq:I}),\allowbreak(\ref{eq:K}),\allowbreak(\ref{eq:Ip}),(\ref{eq:Kp})
appear as eqs.\ (3),(4),\allowbreak(1),\allowbreak(2),(6),(7).  But in
eq.~(4), the factors $e^{\frac12\nu\pi{\rm i}}$, $e^{-\frac12\nu\pi{\rm
    i}}$ should be interchanged.  If reference is made to the definition of
the Ferrers function~${\rm Q}^{\mu}_{\nu}$ in~terms of~${Q}^{\mu}_{\nu}$
\cite[\S\,5.15]{Olver97}, and to the relation between the Barnes and Hobson
definitions of~$Q^{\mu}_{\nu}$, it will be seen that this interchange makes
eq.~(4) equivalent to the simpler identity~(\ref{eq:Y}) of
Theorem~\ref{thm:1}.

\bigskip
\noindent
\textsc{A. Erd\'elyi, W. Magnus, F. Oberhettinger \& F.~Tricomi},
\emph{Tables of Integral Transforms}, Volume~I, McGraw--Hill,
New York, 1954.  Cited here as Ref.~\cite{Erdelyi54}.

\nopagebreak
\medskip
The standard (Hobson) definition of~$Q_\nu^\mu$ is used, and ${\rm
  P}_\nu^\mu,{\rm Q}_\nu^\mu$ are typographically distinct from
${P}_\nu^\mu,{Q}_\nu^\mu$.  On pp.\ 182,187,196,198, versions of 
identities
(\ref{eq:J}),(\ref{eq:Y}),\allowbreak(\ref{eq:I}),\allowbreak(\ref{eq:K}),\allowbreak(\ref{eq:Ip}),(\ref{eq:Kp})
appear as
eqs.\ 4.14(9),\allowbreak4.14(48),\allowbreak4.16(8),\allowbreak4.16(28),\allowbreak4.16(9),\allowbreak4.16(27).

In~4.14(9), $P_\mu^\nu$ should be read as~${\rm P}_\mu^\nu$, making the
identity consistent with~(\ref{eq:J}).  In~4.14(28),
$P_\mu^\nu,P_\mu^{-\nu}$ should be read as~${\rm P}_\mu^\nu,{\rm
  P}_\mu^{-\nu}$.  It follows from the expressions for the second Ferrers
function in \cite[\S\,5.15]{Olver97} that these changes yield consistency
with~(\ref{eq:Y}).

In~4.16(28), the sine factors in numerator and denominator should be
deleted, and $Q_{\mu}^{\nu}$ should be read as $e^{\nu\pi{\rm
    i}}Q_{\mu}^{\nu}$.  In~4.16(9), the sine factors should be deleted, and
$Q_{\nu}^{\mu}$ should be read as $e^{\mu\pi{\rm i}}Q_{\nu}^{\mu}$.  With
these emendations, 4.16(28),4.16(9) will agree with the
identities~(\ref{eq:K}),(\ref{eq:Ip}).  As they stand, they inconsistently
employ the Barnes and not the Hobson definition; they may have been taken
from~\cite{Watson44}.  The corresponding inverse Laplace transforms,
5.13(9) and 5.13(3) on pp.\ 270--271, should be emended similarly.

It has been found that there are many entries in these tables, besides
4.14(9) and 4.14(28), in which Legendre functions are typographically
confused.  A~list is available from the present author.

\bigskip
\noindent
\textsc{G.~E. Roberts \& H. Kaufman},
\emph{Table of Laplace Transforms}, W.~B. Saunders,
Philadelphia, 1966.  Cited here as Ref.~\cite{Roberts66}.

\nopagebreak
\medskip
The standard (Hobson) definition of~$Q_\nu^\mu$ is used, but ${\rm
  P}_\nu^\mu,{\rm Q}_\nu^\mu$ are not visually distinguished from
${P}_\nu^\mu,{Q}_\nu^\mu$.  Versions of identities
(\ref{eq:J}),(\ref{eq:Y}),(\ref{eq:I}),(\ref{eq:K}),(\ref{eq:Ip}),(\ref{eq:Kp})
appear in this table, but its entries 12.2.4 (p.~73) and 13.2.2 (p.~79) are
incorrect, since they incorporate without change the just-noted versions of
(\ref{eq:K}),(\ref{eq:Ip}) found in~\cite{Erdelyi54}.

\bigskip
\noindent
\textsc{M. M. Agrest \& M. S. Maksimov},
\emph{Theory of Incomplete Cylindrical Functions and Their Applications}, Springer-Verlag,
New York/Berlin, 1971.  Cited here as Ref.~\cite{Agrest71}.

\nopagebreak
\medskip
The standard (Hobson) definition of~$Q_\nu^\mu$ is used, but ${\rm
  P}_\nu^\mu,{\rm Q}_\nu^\mu$ are not visually distinguished from
${P}_\nu^\mu,{Q}_\nu^\mu$.  On pp.\ 41--42, versions of identities
(\ref{eq:J}),\allowbreak(\ref{eq:Y}),\allowbreak(\ref{eq:I}),(\ref{eq:K})
appear as eqs.\ (5.3),\allowbreak(5.4),\allowbreak(5.5),(5.6).  In
eq.~(5.5), the prefactor $-(2/\pi)$ should be deleted from the right side.
The domain of validity of eq.\ (5.6) is ${\rm Re}\,a>1$, not ${\rm
  Re}\,a>-1$ as stated; although, the Laplace transform on the left side is
defined if ${\rm Re}\,a>-1$.

Also, in eqs.\ (5.4),(5.5), the representations given for Legendre
functions in~terms of the Gauss hypergeometric function $F={}_2F_1$ are
incorrect: they should be replaced by the representations given
in~\cite{Olver97,Olver2010}. 

In each of eqs.\ (5.8),(5.9), the exponent $\nu+1/2$ should be read
as~$(\nu+1)/2$.  

\bigskip
\noindent
\textsc{I. S. Gradshteyn \& I. M. Ryzhik},
\emph{Table of Integrals, Series, and Products}, 8th ed., Elsevier/Academic
Press, Amsterdam, 2015.  Cited here as Ref.~\cite{Gradshteyn2015}.

\nopagebreak
\medskip
The standard (Hobson) definition of~$Q_\nu^\mu$ is used, but ${\rm
  P}_\nu^\mu,{\rm Q}_\nu^\mu$ are not visually distinguished from
${P}_\nu^\mu,{Q}_\nu^\mu$.  On pp.\ 712--713, in~\S\,6.628, versions of
(\ref{eq:J}),\allowbreak(\ref{eq:Y}),\allowbreak(\ref{eq:I}),\allowbreak(\ref{eq:K}),\allowbreak(\ref{eq:Ip}),(\ref{eq:Kp})
appear as 6.628(1),(2),\allowbreak(4),\allowbreak(5),\allowbreak(6),(7).
All were taken from~\cite{Watson44} without change, so the three containing
the function~$Q_\nu^\mu$ inconsistently employ the Barnes definition,
rather than that of Hobson.  (Equation 6.628(6) has been altered since the
seventh edition, without indication, by correcting a typographical error
that caused it to disagree with~\cite{Watson44}.)

The three requiring emendation are 6.628(2),(5),(6).  Their respective
prefactors,
\begin{displaymath}
  -\,\frac{\sin(\mu\pi)}{\sin[(\mu+\nu)\pi]},\qquad
  \frac{\sin(\mu\pi)}{\sin[(\nu+\mu)\pi]}, \qquad
  \frac{\cos(\nu\pi)}{\sin[(\mu+\nu)\pi]}, 
\end{displaymath}
should be replaced by $-e^{-\nu\pi{\rm i}},e^{-\nu\pi{\rm
    i}},e^{-(\mu-1/2)\pi{\rm i}}$.  

When 6.628(6) is emended as described, it will become identical to the
(correct) integral 6.622(3), so that one or the other will become
redundant.

As was noted above, there is an error in~\cite{Watson44}, in the original
version of eq.\ 6.628(2): $e^{\frac12\nu\pi{\rm i}}$ and
$e^{-\frac12\nu\pi{\rm i}}$ should be interchanged.  If this change too is
made in eq.\ 6.628(2), eq.~6.628(2) will become equivalent to the
identity~(\ref{eq:Y}) listed in Theorem~\ref{thm:1}, which involves the
Ferrers function ${\rm Q}^{-\mu}_\nu$ rather than
${Q}^{\mu}_\nu,{Q}^{-\mu}_\nu$, and is much simpler.


\end{document}